\documentclass[notitlepage, twocolumn, a4paper, prl, nofootinbib, nobibnotes]{revtex4-1}

\usepackage{amsfonts, amsmath, amssymb}
\usepackage[mathcal, mathbf]{euler}
\usepackage[colorlinks=true, dvips, pdfpagelabels]{hyperref}
\usepackage{graphicx}

\def\ulbaddress{Statistical and Plasma Physics, Universit\'e Libre de Bruxelles, Campus Plaine, CP 231, B-1050 Brussels, Belgium}
\def\ucvaddress{Faculty of Physics, University of Craiova, A.I. Cuza 13, 200585 Craiova, Romania}

\begin{document}

	\title{Patterns in the sine map bifurcation diagram}
	\author{C.C. Lalescu}
	\email{clalescu@ulb.ac.be}
	\affiliation{\ulbaddress}
	\affiliation{\ucvaddress}

	\maketitle

	\section{Introduction}
		A onedimensional map $f$ is defined here as any $f:[a,b]\rightarrow[a, b]$.
		The logistic map $\phi_r$
		\begin{equation}
			\phi_r : [0,1] \rightarrow [0, 1], \phi_r (x)= r x (1-x)
		\end{equation}
		is well known for it's chaotic properties.
		As $r$ is varied through different intervals, the logistic map goes from having a fixed point to having stable cycles of different orders $n$ (there is at least an $x$ so that $\phi_r^n(x) = x$, where the exponent refers to composition).
		
\begin{widetext}
			\begin{figure}[h]
				\begin{centering}
					\includegraphics[width=.98\textwidth]{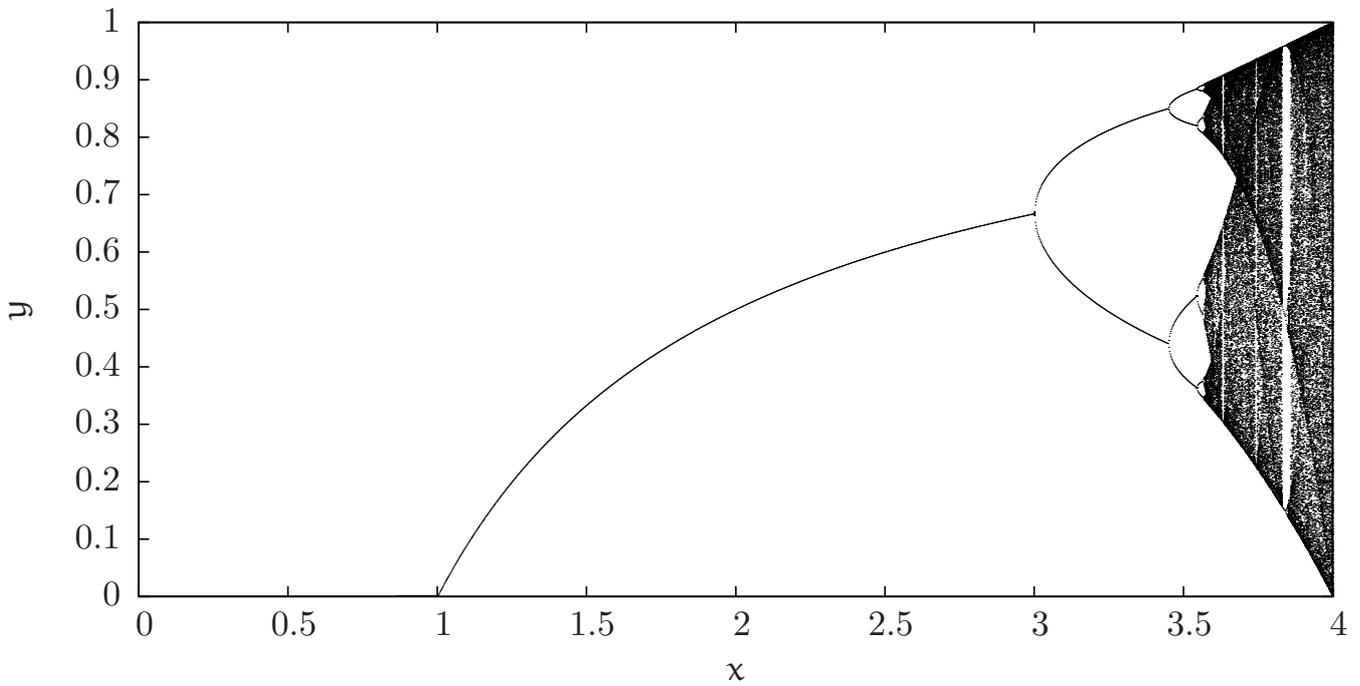}
					\caption{Bifurcation diagram for the logistic map}
					\label{fig:logistic_map}
				\end{centering}
			\end{figure}
		\end{widetext}

		This type of behaviour is generic for a large number of chaotic systems, and it is best viewed in a bifurcation diagram, see figure \ref{fig:logistic_map}.
		The initial fixed point degenerates into a cycle with period 2, and then gradually the period of the stable cycle keeps doubling, up to a point where the map becomes fully chaotic.
		A reference paper discussing the universality of the period doubling cascade is \cite{Feigenbaum_1978}.
		One obvious feature of this bifurcation diagram is that patterns can be observed in it --- discontinuities in the density of points.

	\section{Sine map}
		Consider the map
		\begin{equation}
			\xi_r : \mathbb{R} \rightarrow \mathbb{R}, \xi_r (x)= r \sin(x)
		\end{equation}
		This map generates a bifurcation diagram that is symmetric with respect to both axis of the $(r, x)$ plane.
		\begin{figure*}[t]
			\begin{centering}
				\includegraphics[width=.98\textwidth]{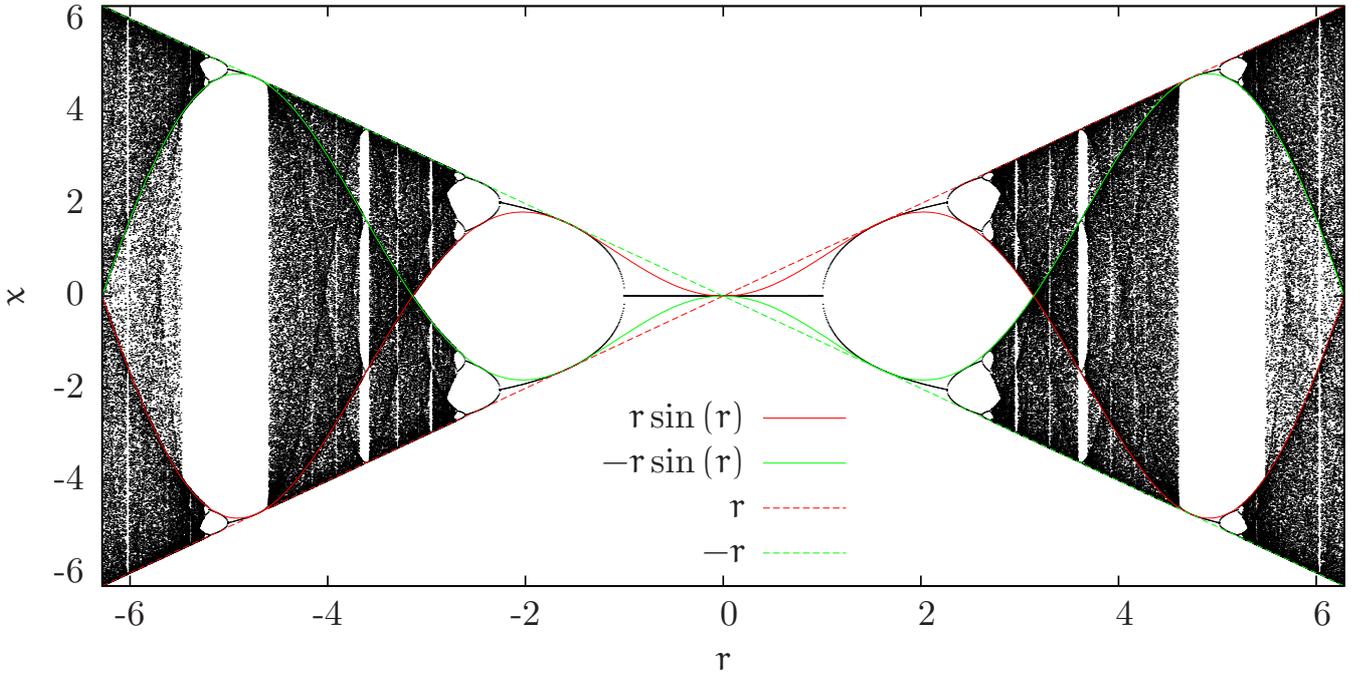}
				\caption{Bifurcation diagram for the sine map, for $r \in [-2\pi,2\pi]$.}
				\label{fig:sine_map}
			\end{centering}
		\end{figure*}
		A twodimensional map can be defined:
		\begin{equation}
			\psi : \mathbb{R}^2 \rightarrow \mathbb{R}^2, \psi (x, y) = (x, x \sin (y))
		\end{equation}
		The bifurcation diagram (viewed as a subset of $\mathbb{R}^2$) is simply the largest subset of $\mathbb{R}^2$ that is invariant to $\psi$:
		\begin{equation}
			\mathcal{B} = \lim_{n \rightarrow \infty} \psi^n (\mathbb{R}^2)
		\end{equation}

	\section{Patterns}
		\begin{figure*}[htb]
			\begin{centering}
				\includegraphics[width=.98\textwidth]{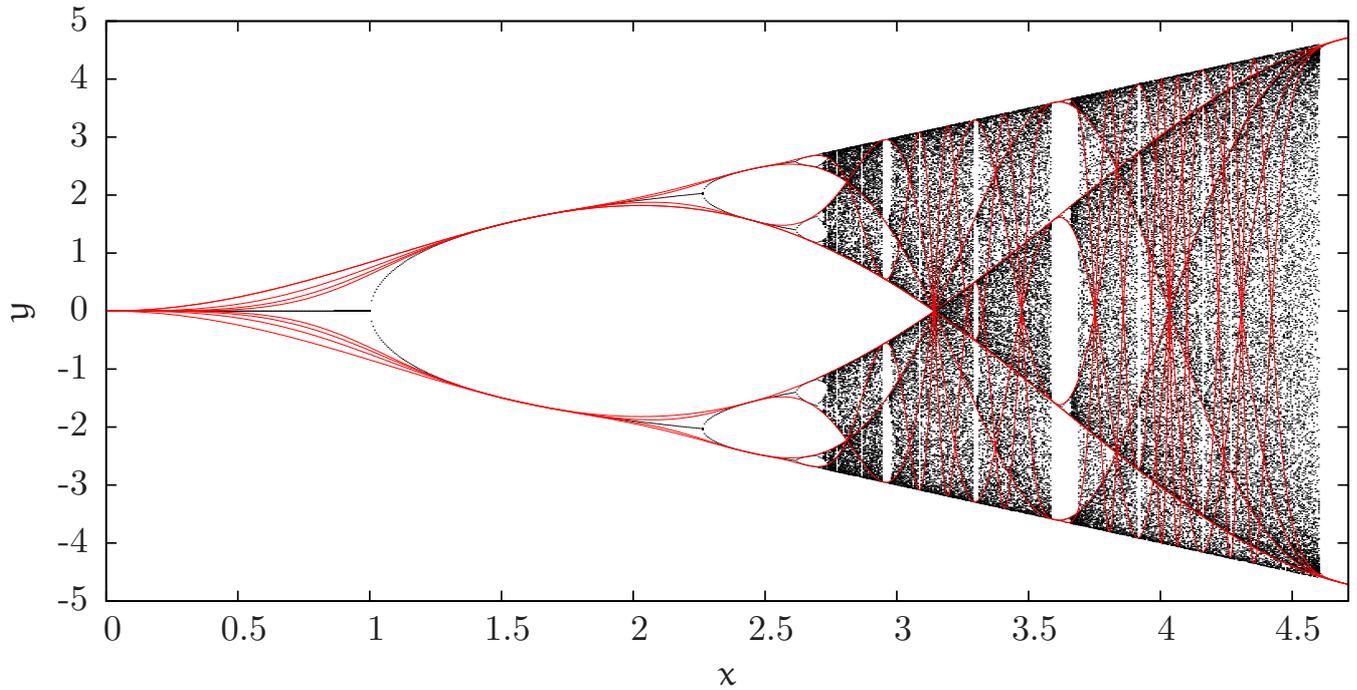}
				\caption{The curve pairs $C^{(n)}$ ($n=1,2,3,4$) superposed over the bifurcation diagram.}
				\label{fig:envelopes}
			\end{centering}
		\end{figure*}
		It can be seen directly in figure \ref{fig:sine_map} that the following two functions appear as patterns:
		\begin{equation}
			\pm r \sin(r), \textrm{ with the graphs given by } \psi(\pm r, r)
		\end{equation}

		More generally, the bifurcation diagrams suggest that the array of curve pairs
		\begin{equation}
			C^{(n)} = \left\{\psi^n(\pm x, x) | x \in \mathbb{R}\right\}
		\end{equation}
		is in fact the succession of all patterns.
		The convention $\psi^0 (x, y) \equiv (x, y)$ is adopted here, so $C^{(0)} = \{(x, x) , (-x, x) | x \in \mathbb{R}\}$.
		Figure \ref{fig:envelopes} further suggests that $C^{(n)} \cap C^{(m)} \subset \mathcal{B}$, for all $n\neq m \in \mathbb{N}$.

	\section{Exact results}
		Let $(a, b) \in C^{(n)} \cap C^{(m)}$.
		It is convenient to introduce the notation $\psi^n_\pm (x, y) \equiv \psi^n(\pm x, y)$.
		This means that
		\begin{equation}
			(a, b) \in \{ \psi^n_\pm(a, a), \psi^m_\pm (a, a)\}
		\end{equation}
		Consider the case $(a, b) = \psi^n_+(a, a) = \psi^m_+(a, a)$.
		Assume $m>n$, and $m - n = p$.
		This means that
		\begin{equation}
			\psi^n_+(a, a) = \psi^p_+( \psi^n_+(a, a))
		\end{equation}
		so, by definition, $\psi^n_+(a, a) \equiv (a, b)$ is part of a cycle with period $p$, which means it is in $\mathcal{B}$.
		The alternative is that $\psi^n_+(a,a) = \psi^m_-(a,a)$.
		This means that
		\begin{align}
			\psi^n_+(a, a) &= \psi^p_+( \psi^n_-(a, a)) \\
			\psi^p_-(\psi^n_+(a, a)) &= \psi^p_-(\psi^p_+( \psi^n_-(a, a))) \\
			\psi^p_+(\psi^n_-(a, a)) &= \psi^p_+(\psi^p_+( \psi^n_+(a, a))) \\
			\psi^n_+(a, a) &= \psi^{2p}_+( \psi^n_+(a, a))
		\end{align}
		so, again by definition, $(a, b)$ is part of a cycle with period $2p$, so $(a, b) \in \mathcal{B}$.

		Concerning $\lim_{n \rightarrow \infty} C^{(n)}$, it is sufficient to remark that $C^{(n)} = \psi^n (C^{(0)}$, so the limit $C^{(\infty)}$ is a subset of $\mathcal{B}$.

	\section{Generalization}
		Consider a general mapping
		\begin{equation}
			\theta : \mathbb{R}^2 \rightarrow \mathbb{R}^2, \theta(x, y) = (x, x f(y))
		\end{equation}
		For $f:\mathbb{R} \rightarrow [-1, 1]$ (with $f(\mathbb{R}) = [-1, 1]$), the invariant set of $\theta$ is well defined as
		\begin{equation}
			\mathcal{D} \equiv \lim_{n\rightarrow \infty} \theta^n\left(\mathbb{R}^2\right)
		\end{equation}
		and the reasoning from the previous section can be applied.
		Call $\theta_\pm (x, y) \equiv \theta ( \pm x, y)$ and
		\begin{equation}
			D^{(n)} \equiv \left\{ \theta^n_\pm (x, x) \big| x \in \mathbb{R} \right\}
		\end{equation}
		Then for any $n\neq m \in \mathbb{N}$, $D^{(n)} \cap D^{(m)} \subset \mathcal{D}$ and the limit $D^{(\infty)}$ is a subset of $\mathcal{D}$.

	\begin{widetext}
		\section{Notes}
			The author would like to thank Gy. Steinbrecher and D. Constantinescu from the University of Craiova for their valuable comments.
			
			The author is unaware if the ideas presented are new; the webpage \cite{Evgeny_Demidov} does discuss a family of curves that appears to have the properties mentioned in this work.
			The figures in this article were generated with a simple Python script\footnote{\url{http://chichi.lalescu.ro/files/bc.py}}.

	\end{widetext}

\end{document}